                    \def\version{20 October, 2026}                       %

 \documentclass[reqno,11pt]{amsart}
 \usepackage{amsmath, amsthm, a4, latexsym, amssymb, mathtools, enumitem, xcolor, color, tikz-cd,todonotes}
\usepackage[unicode]{hyperref}
\usepackage{srcltx}

\def\@rmrk#1#2{\refstepcounter
    {#1}\@ifnextchar[{\@yrmrk{#1}{#2}}{\@xrmrk{#1}{#2}}}

\makeatletter\@addtoreset{equation}{section}\makeatother

 \newfont{\bfit}{cmbxti10 scaled 1200}

\renewcommand{\d}{{\rm d}}
 \newcommand{\e}{{\rm e} }

 \newcommand{\eps}{\varepsilon}
 \newcommand{\supp}{{\rm supp}}
  \newcommand{\Aut}{{\rm Aut}}

 \newcommand{\R}{\mathbb{R}}
\newcommand{\N}{\mathbb{N}}

\newcommand{\C}{\mathbb{C}}

\newcommand{\E}{\mathbf{E}}
\renewcommand{\P}{\mathbf{P}}
\newcommand{\PP}{\mathbb{P}}

 \def\1{{\mathchoice {1\mskip-4mu\mathrm l} 
{1\mskip-4mu\mathrm l}
{1\mskip-4.5mu\mathrm l} {1\mskip-5mu\mathrm l}}}

 \newcommand{\skrib}{{\mathcal B}}
 \newcommand{\skric}{{\mathcal C}}

\newcommand{\heap}[2]{\genfrac{}{}{0pt}{}{#1}{#2}}

\renewcommand{\subsection}{\secdef \subsct\sbsect}
\newcommand{\subsct}[2][default]{\refstepcounter{subsection}
\vspace{0.15cm}
{\flushleft\bf \arabic{section}.\arabic{subsection}~\bf #1  }
\nopagebreak\nopagebreak}
\newcommand{\sbsect}[1]{\vspace{0.1cm}\noindent
{\bf #1}\vspace{0.1cm}}

{\nopagebreak {\hfill\rule{2mm}{2mm}}\\ }

\newtheorem{theorem}{Theorem}[section]
\newtheorem{lemma}[theorem]{Lemma}
\newtheorem{cor}[theorem]{Corollary}
\newtheorem{prop}[theorem]{Proposition}

\newtheorem{theorem*}{Theorem}
\newcommand{\eproof}{\hfill \qed \vspace*{5mm}}

\newtheoremstyle{thm}{1.5ex}{1.5ex}{\itshape\rmfamily}{}
{\bfseries\rmfamily}{}{2ex}{}

\newtheoremstyle{rem}{1.3ex}{1.3ex}{\rmfamily}{}
{\itshape\rmfamily}{}{1.5ex}{}
\theoremstyle{rem}
\newtheorem{remark}{{\slshape\sffamily Remark}}[]

\refstepcounter{subsubsection}

\def\thebibliography#1{\section*{References}
  \list%
  {\arabic{enumi}.}
    {\settowidth\labelwidth{[#1]}\leftmargin\labelwidth
    \advance\leftmargin\labelsep
    \parsep0pt\itemsep0pt
    \usecounter{enumi}}
    \def\newblock{\hskip .11em plus .33em minus .07em}
    \sloppy                   
    \sfcode`\.=1000\relax}

\begin{document}
\title[Coarse embeddability, $L^1$-compression and Percolations on General Graphs]
{\large Coarse embeddability, $L^1$-compression and Percolations on General Graphs}
\author[Chiranjib Mukherjee and Konstantin Recke]{}
\maketitle
\thispagestyle{empty}
\vspace{-0.5cm}

\centerline{By {\sc Chiranjib Mukherjee\footnote{Universit\"at M\"unster, Einsteinstrasse 62, M\"unster 48149, Germany, {\tt chiranjib.mukherjee@uni-muenster.de}}}
and {\sc Konstantin Recke\footnote{Universit\"at M\"unster, Einsteinstrasse 62, M\"unster 48149, Germany {\tt konstantin.recke@uni-muenster.de}}}}
\renewcommand{\thefootnote}{}
\footnote{Corresponding author: Chiranjib Mukherjee, {\tt chiranjib.mukherjee@uni-muenster.de}}

\vspace{-0.5cm}
\centerline{\textit{Universit\"at M\"unster}}
\vspace{0.2cm}

\begin{center}
\version
\end{center}

\begin{quote}{\small {\bf Abstract:} We show that a locally finite, connected graph has a coarse embedding into a Hilbert space if and only if there exist bond percolations with arbitrarily large marginals and two-point function vanishing at infinity. We further show that the decay of the two-point function is stretched exponential with stretching exponent $\alpha\in[0,1]$ if and only if the $L^1$-compression exponent of the graph is at least $\alpha$, leading to a probabilistic characterization of this exponent. 
These results are new even in the particular setting of Cayley graphs of finitely generated groups. 
The proofs build on new probabilistic methods introduced recently by the authors to study group-invariant percolation on Cayley graphs \cite{MR22,MR23}, which are now extended to the general, non-symmetric situation of graphs to study their coarse embeddability and $L^1$-compression exponents.} 
\end{quote}

\section{Main results} \label{sec intro}

\subsection{Background.}\label{section-background}\noindent Consider a locally finite, connected, infinite graph $G=(V,E)$. The incentive to understand the {\em large scale geometry} of $G$ leads to the two following important notions of measuring to which extent its geometry is comparable to that of a Hilbert space, resp.~to that of an $L^1$-space:

\begin{itemize}
    \item[$\bullet$] A map $f\colon X\to Y$ between metric spaces is a {\bf coarse embedding}, if there exist non-decreasing maps $\rho_1,\rho_2\colon[0,\infty)\to[0,\infty)$ such that $\rho_1(t)\to\infty$ as $t\to\infty$ and
    \begin{equation} \label{equ-CoarseEmbedding}
    \rho_1(d_X(u,v))\leq d_Y(f(u),f(v)) \leq \rho_2(d_X(u,v))
    \end{equation}
    for all $u,v\in X$. We are interested in the existence of a coarse embedding from $X=G$ into the target space $Y=H$ for some Hilbert space $H$.
    
    \vspace{1mm}
    
    \item[$\bullet$] The {\bf $L^1$-compression exponent} $\alpha_1^*(X)$ of a metric space $X$ is the supremum over all those $\alpha\in[0,1]$ such that there exist a measure space $(Y,\nu)$ and a Lipschitz function $f\colon X\to L^1(\nu)$ satisfying 
    \begin{equation}\label{equ-DefL1Compression}
    \Vert f(u)-f(v) \Vert_1 \geq c d_X(u,v)^\alpha
    \end{equation}
    for every $u,v\in X$ and a uniform constant $c>0$, see \cite{GK04,NP11}. We are interested in $\alpha_1^*(G)$, which provides an elegant measurement of the non-bi-Lipschitz embeddability of the graph into $L^1$-spaces.
\end{itemize}
These concepts fit into a broader program of embedding complicated structures into target geometries with controlled distortion. For instance, embeddings of finite metric spaces into (finite dimensional) Hilbert spaces and $\ell^p$-spaces have been studied extensively in combinatorics. The above are extensions to the setting of infinite graphs, where a coarse embedding may be understood as a minimal assumption on controlling the distortion, while the compression exponent provides a quantitative measurement of the failure of bi-Lipschitz embededdability. 

Besides their intrinsic interest, these concepts have notable applications to group theory and probability theory. The property of being coarsely embeddable into a Hilbert space is perhaps most known for Yu's \cite{Y00} seminal work showing that it implies the coarse Baum-Connes and Novikov conjectures. The $L^1$ and similarly defined $L^p$ compression exponents are well-known quasi-isometry invariants introduced by Guentner and Kaminker \cite{GK04}. They proved for instance that an $L^2$-compression exponent strictly larger than $1/2$ implies Yu’s property~A. Note that there are groups which coarsely embed into Hilbert space with $L^2$-compression exponent equal to zero~\cite{ADS09}. The case $p=1$ is perhaps the most natural beyond $p=2$ and still quite rich -- notably, there are amenable groups with $\alpha^*_1=0$ by the work of Austin~\cite{A11}. The $L^1$-target spaces also have particularly appealing features: they are related to cut-metrics well-known in combinatorics and the measure definite kernels introduced by Robertson and Steeger~\cite{RS98} (to study infinite measure preserving actions of groups with Kazhdan’s property~(T)). Finally, another main motivation for us to study these concepts using probabilistic techniques came from the significant activity around random walk approaches to $L^p$-compression exponents, see e.g.~\cite{NP11,BZ21} and references therein.

In this paper, we develop a novel method using percolations on $G$ leading to the first probabilistic characterizations of coarse embeddability into a Hilbert space and of the $L^1$-compression exponent.
More precisely, we characterize these properties by the existence of random subgraphs obtained by randomly deleting or retaining edges (resp.~vertices) -- that is by  considering {\em bond (resp.~site) percolations} -- which simultaneously satisfy two competing properties: 
\begin{itemize}
    \item[(i)] containing every individual edge (resp.~site) with high probability, {\em and} 
    
    \vspace{1mm}
    
    \item[(ii)] exhibiting decay of connection probabilities (i.e.~decay of the two-point function), the decay rate being quantified explicitly.
    \end{itemize}

To precisely state our results, recall that a {\bf general bond percolation} on a graph $G=(V,E)$ is the distribution of a random subgraph obtained by keeping or deleting edges. We say that a general bond percolation $\P$ has {\bf marginals larger than} $p\in[0,1]$, if 
\begin{equation} \label{equ-LargeMarginals}
    \inf_{e\in E} \, \P\big[e\in E] \geq p.
\end{equation}
The {\bf two-point function}
\begin{equation} \label{def-TPF}
\tau \colon V \times V \to [0,1] \, , \, \tau(u,v):=\P\big[u\leftrightarrow v\big]
\end{equation}
will be our main object of interest. Its decay or non-decay at infinity measures the {\em connectivity} of the percolation under consideration. Note the competitive relationship between having large marginals (i.e.~larger than some $p$ which is close to $1$) and simultaneously exhibiting two-point function decay. This leads to the natural problem of determining under which conditions, and for which kinds of models, two-point function decay may hold, or may fail, in the presence of large marginals.

In the setting of Cayley graphs of groups, several important geometric properties have recently been shown to admit characterizations \cite{MR22,MR23} through the existence or non-existence of {\em invariant} percolations (i.e.~invariant under the natural action of the group on its Cayley graph) with large marginals and two-point function decay. More precisely, the Haagerup property and Kazhdan's property~(T) have been characterized through percolation by the present authors \cite{MR23} -- we refer to Section \ref{section-ClosingRemarks} for details. For now, let us record that these results assert that on highly symmetric graphs, connectivity of invariant percolations is closely related to geometric properties of the symmetry group. In this paper, we will provide a surprising extension of the method developed by the authors in \cite{MR23} to obtain results in the -- in a sense completely orthogonal -- situation of coarse embeddability and compression exponents of {\em non-symmetric} graphs;~see Theorem \ref{theorem-CE} and Theorem \ref{theorem-L1NE} for precise statements. Their consequences (Corollary \ref{cor-CE-groups} and Corollary \ref{cor-L1NE-groups}) are new also in the particular setting of Cayley graphs of finitely generated groups. We also refer to Remark~\ref{rem-Lp} and to Section~\ref{section-ClosingRemarks}~(B)for applications of these results to $L^p$-compression exponents for $p>1$. These results in particular characterize graphs for which large marginals and (stretched exponential) decay of two-point functions can be reconciled, and provide examples when this necessarily fails for marginals above a critical threshold. 

A main appeal of the current method is, therefore, the great generality in which it applies while relying on fairly simple, yet conceptually rich and intuitive probabilistic ideas. Given its simplicity, it is likely that this approach will facilitate a wider range of new applications to the coarse geometric properties of graphs using these probabilistic concepts. Let us now turn to the formal statements of our results.

\subsection{Results about general infinite graphs.} As described above, we study two-point function decay of general bond percolations on general infinite graphs. These of course form a rather wild set and it is a priori far from clear whether the corresponding two-point function decay encodes {\em significant information at all} (this is elaborated in Remark \ref{rem-GenPerc}). Our first main result shows that,  in fact, it does: 

\begin{theorem}[Percolation and coarse embeddability into a Hilbert space for graphs] \label{theorem-CE}
Let $G=(V,E)$ be a locally finite, connected graph with $d$ denoting the graph metric. Then $G$ admits a coarse embedding into a Hilbert space (recall~\eqref{equ-CoarseEmbedding}) if and only if for every $p<1$, there exists a general bond percolation $\P$ on $G$ with $\P\big[e\in\omega\big]>p$ for every $e\in E$ and such that the two-point function vanishes at infinity, i.e.
\begin{equation} \label{equ-TPFvanishing}
\lim_{r\to\infty} \sup_{\heap{u, v \in V}{d(u,v)>r}} \,\, \P\big[u\leftrightarrow v\big]  = 0.
\end{equation}
$\mathrm{(See \,\,Theorem \,\,\ref{theorem-Coarse-Embeddability}).}$
\end{theorem}

In particular, we obtain the following consequence in the invariant setting. Here a finitely generated group {\bf embeds coarsely into a Hilbert space}, if there exists a coarse embedding into a Hilbert space in the sense of \eqref{equ-CoarseEmbedding} for the group equipped with some word metric.

\begin{cor}[Percolation and coarse embeddability into a Hilbert space for groups] \label{cor-CE-groups} Let $\Gamma$ be a finitely generated group. Then $\Gamma$ embeds coarsely into a Hilbert space if and only if some, equivalently every, Cayley graph $G=(V,E)$ has the property that for every $p<1$, there exists a general bond percolation $\P$ with $\P\big[e\in\omega\big]>p$ for every $e\in E$ and such that the two-point function vanishes at infinity.
\end{cor}

The following remarks will put Theorem \ref{theorem-CE} and Corollary \ref{cor-CE-groups} into context.

\begin{remark}(Invariance vs.~non-invariance) \label{rem-GenPerc} We emphasize that there is no assumption of invariance in Theorem \ref{theorem-CE}, which is quite exceptional. To quote from a beautiful survey by H{\"a}ggstr{\"o}m and Jonasson~\cite{HJ06}, the assumption of group-invariance is ``extremely natural'' and ``almost universally employed'' in percolation theory, and ``to work with arbitrary probability measures on $\{0,1\}^E$ would, however, be to take things a bit too far, as not much of interest can be said in such a general setting". This intuition is underscored by a plethora of results in the invariant setting, see \cite{LP16,P22} (see also Remark~\ref{rem-A} for the only exception the authors are aware of). From this point of view, it is quite surprising that two-point function decay of general percolations characterizes coarse embeddability into a Hilbert space, i.e.~the natural notion of looking like a Hilbert space on large scales. This is particularly noteworthy in Corollary \ref{cor-CE-groups} because there is a canonical group action available. In fact, the statement with  $\Gamma$-invariant bond percolations characterizes the Haagerup property of~$\Gamma$ (cf.~Section~\ref{section-ClosingRemarks}~(A)). 
In the setting of groups, this relationship between coarse embeddability and the Haagerup property (which can be interpreted as the existence of an equivariant coarse embedding into~$L^2$; see~Section~\ref{section-ClosingRemarks}~(A)) recovers similar observations in coarse geometry, see e.g.\ \cite{BO08,W09}.
\nopagebreak {\hfill\rule{2mm}{2mm}}
\end{remark}

\begin{remark}[Probabilistic interpretation of coarse embeddability] \label{rem-CE} To the best of the authors' knowledge, Theorem \ref{theorem-CE} is the first probabilistic characterization of coarse embeddability into a Hilbert space. Its geometric meaning may be understood by thinking of the graph as a network and considering the task of designing a strategy of randomly removing edges in order to disconnect the graph. Consider the strategy successful if vertices, which are far apart, are disconnected with high probability, i.e.~with probability tending to one as the distance tends to infinity. As a constraint, associate a cost to the removal of edges so that it is only allowed to remove each edge with probability at most $1-p$. Theorem~\ref{theorem-CE} shows that there exists a successful strategy at arbitrarily small costs if and only if the graph looks like a Hilbert space on large scales. \nopagebreak {\hfill\rule{2mm}{2mm}}
\end{remark}

\begin{remark}[Difference to property A]\label{rem-A} For a {\it bounded degree} graph, coarse embeddability into a Hilbert space follows from Yu's property A, introduced in \cite{Y00}. Following \cite{W09}, we say that $G=(V,E)$ has {\bf property A} if for all $R,\eps>0$, there exists a family $\{A_v\}_{v\in V}$ of non-empty finite subsets of $V\times\N$ such that $d(u,v)\leq R$ implies $|A_u \Delta A_v|/|A_u\cap A_v|<\eps$, and such that there exists $S\geq0$ such that $(u,n)\in A_v$ implies $d(v,u)\leq S$. 
It has been shown in a series of works (e.g.~\cite[Section 2]{E21b}) that a connected graph $G=(V,E)$ with bounded degrees has property A if and only if for every $\eps>0$ there exists $K\geq1$ and a probability measure $\nu$ on the set of $K$-separators (i.e.~subsets $A\subset V$ such that all components of the induced graph on $V\setminus A$ are of size at most $K$) such that for every $v\in V$, the measure of the set of $K$-separators containing $v$ is at most~$\eps$.
By considering edges in the complement of a random $K$-separator, we may state an equivalent condition in the present paper's terminology as follows: for every $p<1$, there exists $K\geq1$ and a general bond percolation $\P$ with $\P\big[e\in\omega\big]>p$ for every $e\in E$ and such that all clusters are of size at most $K$ (this condition may equivalently be replaced by requiring the two-point function to be supported on finite tubes of diameter $K$, i.e.~$\P\big[u\leftrightarrow v\big]=0$ whenever $d(u,v)>K$). Although Theorem \ref{theorem-CE} resembles this in spirit, the methods used for property~A are quite different from the present paper and also do not seem to apply in our setting -- in particular, Theorem~\ref{theorem-CE} does not assume bounded degrees. 
\nopagebreak {\hfill\rule{2mm}{2mm}}
\end{remark}

\begin{remark}[Amenability] \label{rem-Amenability}
To elaborate on the previous point, 
property~(A) is similar to amenability in the sense that it is characterized intrinsically by the existence of sets with small boundaries, which is not the case for coarse embeddability into a Hilbert space. Amenability for quasi-transitive graphs is a fundamental property, which has been characterized through the existence of group-invariant percolations with large marginals and only finite clusters in the seminal work of Benjamini, Lyons, Peres, and Schramm \cite{BLPS99}, see also \cite{MR22} for a closely related characterization emphasizing the role of the two-point function. These results rely on a powerful tool in percolation known as the {\em mass-transport principle}~\cite{BLPS99}. We do not use the mass-transport principle to prove Theorem \ref{theorem-CE} and point out that this principle does not hold at the present level of generality.
\nopagebreak {\hfill\rule{2mm}{2mm}}
\end{remark}

Let us now turn to the main result concerning {\em compression exponents}, which provide a way of measuring non-bi-Lipschitz embeddability introduced in \cite{GK04}: recall that the {\bf $L^1$-compression exponent} $\alpha_1^*(X)$ of a metric space $X$ is the supremum over those $\alpha\in[0,1]$ such that there exist a measure space $(Y,\nu)$ and a Lipschitz function $f\colon X\to L^1(\nu)$ satisfying 
\begin{equation*}
\Vert f(u)-f(v) \Vert \geq c d_X(u,v)^\alpha
\end{equation*}
for every $u,v\in X$ and a uniform constant $c>0$, cf.~\cite{GK04,NP11}. Our second main result characterizes $\alpha_1^*$ of general graphs as well as finitely generated groups through percolation:

\begin{theorem}[Percolation and the $L^1$-compression exponent for graphs]\label{theorem-L1NE} 
Let $G=(V,E)$ be a locally finite, connected graph. Then the $L^1$-compression exponent $\alpha_1^*(G)$ is the supremum over those $\alpha\in[0,1]$ for which there exists $C>0$ such that for every $p<1$, there exists a general bond percolation $\P$ on $G$ with $\P\big[e\in\omega\big]>p$ for every $e\in E$ and which satisfies the stretched exponential decay
\begin{equation}\label{equ-TPFstretchedexp}
\e^{-\beta d(u,v)} \leq \P\big[u \leftrightarrow v \big] \leq \e^{-\gamma d(u,v)^\alpha} 
\end{equation}
for every $u,v\in V$ and some $\beta,\gamma>0$ with $\beta/\gamma\leq C$.

$\mathrm{(See \,\, Theorem \,\, \ref{theorem-L1-Compression-Exponent}).}$
\end{theorem}

In fact,  we will provide a more general result for arbitrary compression functions, i.e.~lower bounds in~\eqref{equ-DefL1Compression} which are of the form $c\rho(d(u,v))$ for choices of $\rho$ possibly different from $\rho(t)=t^\alpha$ (see Theorem~\ref{theorem-L1-Compression-Functions} for details). We refer to Remark~\ref{rem-Lp} for applications of this result to $L^p$-compression exponents with $p>1$ and to Section~\ref{section-ClosingRemarks}~(B) for applications to equivariant $L^p$-compression exponents with $p\ge1$. As an important corollary, we characterize the $L^1$-compression exponent of a finitely generated group, which is defined as in \eqref{equ-DefL1Compression} for the group equipped with some word metric, cf.~\cite{NP11}.

\begin{cor}[Percolation and the $L^1$-compression exponent for groups] \label{cor-L1NE-groups}
Let $\Gamma$ be a finitely generated group and let $G=(V,E)$ be some Cayley graph. Then $\alpha_1^*(\Gamma)$ is the supremum over those $\alpha\in[0,1]$ for which there exists $C>0$ such that for every $p<1$, there exists a general bond percolation $\P$ with $\P\big[e\in\omega\big]>p$ for every $e\in E$ and which satisfies the stretched exponential decay
\begin{equation*}
\e^{-\beta d(u,v)} \leq \P\big[u \leftrightarrow v \big] \leq \e^{-\gamma d(u,v)^\alpha} 
\end{equation*}
for every $u,v\in V$ and some $\beta,\gamma>0$ with $\beta/\gamma\leq C$.
\end{cor}

Let us conclude with two remarks, which provide the necessary context.

\begin{remark}[Comparison with equivariant compression]\label{remark comparison}
In contrast to the invariant situation considered in \cite{MR23} (see~Section~\ref{section-ClosingRemarks}~(B) for details), Theorem \ref{theorem-L1NE} shows that the $L^1$-compression exponent is {\em characterized} by the existence of bond percolations with large marginals and sufficiently strong control of the two-point function decay. This significant strengthening is enabled by a clearer understanding of the kernel-theoretic technicalities which show up in the proof, see Section \ref{section-Kernels} and Remark \ref{rem-L1andL1NE}. \nopagebreak {\hfill\rule{2mm}{2mm}}
\end{remark}

\begin{remark}[Probabilistic interpretation of the $L^1$-compression exponent]\label{rem-L1ProbInt}  
Theorem \ref{theorem-L1NE} characterizes the $L^1$-compression exponent of locally finite, connected graphs as well as finitely generated groups by a probabilistic robustness property of the graph similarly to Remark \ref{rem-CE}. More precisely, it shows that a graph does not admit an embedding into an $L^1$-space with good compression if and only if any strategy of disconnecting the graph fails to deliver stretched exponential decay in the sense of~\eqref{equ-TPFstretchedexp} for marginals above some non-trivial threshold.
Let us also point out that there exist powerful interactions between the theory of compression exponents and random walks, see e.g.~\cite{BZ21,NP08,NP11} and the references there -- our results open a path towards establishing further links with percolation. \nopagebreak {\hfill\rule{2mm}{2mm}}
\end{remark}

This concludes our investigation of the connections between two-point function decay of percolations and geometric properties of general graphs. Let us also include the brief comment, that the preceding two results perhaps suggest, to a reader familiar with one of the results in \cite{MR23} (see Remark (A) and Remark (E) in Section \ref{section-ClosingRemarks} for relevant background) that a non-invariant version of Kazhdan's property~(T) can also be formulated 
-- we show in Proposition~\ref{prop-GeneralDecay} 
that this is not the case. Indeed, there it is shown that every locally finite, connected, infinite graph admits percolations (satisfying the FKG-inequality) with large marginals and two-point function decay {\em in some direction}.

\subsection{Ideas of the proofs.} Let us conclude this section by highlighting the two main difficulties one needs to overcome to obtain the present results. First, one would expect a statement about the behavior of all percolation models, and in particular a statement about their two-point function decay, to be essentially vacant in the sense that it would not bear significant information about interesting geometric features of the graph. Secondly, if one were to set about understanding embeddings into Hilbert spaces or $L^1$-spaces through percolations, the natural attempt would be to construct percolation models capturing the geometry {\em intrinsically}, i.e.~by removing the boundaries of certain sets randomly. This approach, which works very well for amenability due to the existence of F{\o}lner sets, can not be implemented for obtaining embedding properties, cf.~Remark \ref{rem-Amenability}.

As we will see, the way to address these difficulties is to use a novel machinery relating {\em measure definite kernels} and percolations via monotone couplings of Poisson point processes (in a non-invariant setting). This machinery makes it possible to obtain results about relevant geometric features despite the absence of a group action, which is perhaps the first distinguishing feature of our approach. Moreover, this machinery is especially robust in the present setting of locally finite, connected graphs, which we emphasize is perhaps as general as one could hope for.

The groundwork for our approach is laid out in Section \ref{section-Preliminaries}, in which, following a brief review of the relevant theory of kernels, we provide a key technical tool which links measure definite kernels and percolation. Namely, we show in Lemma~\ref{lemma-TPF-MD} that every two-point function $\tau$ induces the measure definite kernel $1-\tau$. This property leads to a new ``duality'' between two-point functions and measure definite kernels. More precisely, building on the strategy developed by the authors for invariant percolation on Cayley graphs \cite{MR23}, we construct for every measure definite kernel a monotonically coupled family $(\P_t)_{t>0}$ of general percolations  satisfying the FKG-inequality, having marginals interpolating between $0$~and~$1$ as $t\to0$ and having explicit control of the two-point function decay at infinity in terms of the growth of the kernel, see Corollary \ref{cor-ConstructionL1}. To define these percolations, we use a (non-invariant) Poisson point process on an infinite measure space associated to every measure definite kernel (see Lemma~\ref{lemma-CSV}). The usefulness and potential of this machinery is underscored by the fact that it is the key ingredient in the proofs of both main results. Concretely, in the proof of Theorem~\ref{theorem-CE} we apply it to derive from a characterization of coarse embeddability through kernels our probabilistic characterization. Here, the construction of the desired percolations is the more challenging direction of the characterization, which is taken care of by the aforementioned construction. On the other hand, the proof of Theorem \ref{theorem-L1NE} is where this machinery is most useful and improves significantly analogous results for invariant percolation on Cayley graphs (see Remark \ref{rem-L1andL1NE}). Here, good $L^1$-compression can be used to prove existence of percolations with stretched exponential two-point function decay through our construction, while the existence of such percolations implies the existence of measure definite kernels with large growth (which implies good $L^1$-compression) through the aforementioned duality.

\subsection{Organization.} The rest of this paper is organized as follows: Section \ref{section-Preliminaries} introduces terminology and notation regarding graphs and percolation and provides the background for a general theory of kernels. In Section \ref{sec MD2PWall} we will establish the aforementioned duality between two-point functions of general percolations and measure definite kernels, while in Section \ref{section-ConstructionWalls} the measure theoretic construction about wall structures on general graphs will be provided. Section \ref{section-CEandL1} contains the construction of one-parameter family of percolations using non-invariant Poisson point processes (see Section \ref{section-ConstructionPerc} for the definition and properties of the constructed percolations) as well as the proofs of our main results (see Section \ref{section-proofs}). In Section \ref{section-ClosingRemarks}, we provide a brief outlook and questions raised by our results.

\section{Preliminaries} \label{section-Preliminaries}

This section introduces the key notions from graph theory, percolation and the theory of kernels needed in the sequel. For more extensive treatments of these topics, we refer to \cite{BdlHV08,BO08,LP16,P22,W09} (see also \cite{MR23} for a more detailed discussion of the invariant setting and \cite{L13} for a different probabilistic treatment of kernels).

\subsection{Graph theoretic terminology and percolation.} \label{section-GraphTerminology} Consider a graph $G=(V,E)$ with vertex set $V$ and symmetric edge set $E\subset V\times V$. An edge between two vertices $u$ and $v$ is denoted by $[u,v]$. In this case, we say that $u$ and $v$ are {\bf adjacent} or {\bf neighbors} and write $u\sim v$. The {\bf degree} $\deg(v)=\deg_{G}(v)$ of a vertex $v$ is the number of vertices adjacent to $v$. If $H$ is a subgraph of $G$, $\deg_H(v)$ denotes the number of adjacent vertices in $H$. The graph is {\bf locally finite}, if every vertex has finite degree. It has {\bf bounded degree} if $\sup_{v\in V} \deg(v)<\infty$. The graph is {\bf connected}, if every two vertices are joined by a path of edges. For a connected graph, the {\bf graph distance} between $u,v\in V$ is defined to be the length of a shortest path connecting $u$ and $v$ and is denoted by $d(u,v)$.

 A {\bf general bond percolation} on a graph $G=(V,E)$ is a probability measure on subsets of $E$, i.e.~on $\{0,1\}^E$ equipped with the Borel $\sigma$-algebra. The percolation configuration is typically denoted by $\omega$ and we identify $\omega$ with the corresponding induced subgraph without further mention. A {\bf cluster} of $\omega$ is a connected component of $\omega$ and we write $C(v)=C_\omega(v)$ for the cluster containing $v\in V$. We write $\deg_\omega(v)$ for the degree of $v\in V$ as a vertex in $\omega$. The {\bf two-point function} 
\begin{equation}
\tau \colon V \times V \to [0,1] \, , \, \tau(u,v):=\P[u\leftrightarrow v],
\end{equation}
is the probability that $u$ and $v$ are in the same cluster of the percolation configuration. If $\Gamma\subset\Aut(G)$ is a subgroup, we say that $\P$ is {\bf$\Gamma$-invariant} if it is invariant under the induced action of $\Gamma$ on $\{0,1\}^E$. The most prominent (invariant) model is {\bf Bernoulli bond percolation} $\PP_p$ with parameter $p\in[0,1]$, in which each edge is deleted independently with probability $1-p$.
Bernoulli percolation is well-known to satisfy the FKG-inequality, which we now introduce in general.

An event $A\in\skrib(\{0,1\}^E)$ is {\bf increasing}, if $\omega\in A$ and $\omega\leq\omega'$ implies $\omega'\in A$. We say that $\P$ satisfies the {\bf FKG-inequality} (or has {\bf positive associations}), if for any two increasing events $A$ and $B$
\begin{equation} \label{FKG-inequality}
\P(A\cap B) \geq \P(A) \P(B).
\end{equation}
Similarly, a measurable funcion $f\colon\{0,1\}^E\to\R$ is {\bf increasing} if $\omega\leq\omega'$ implies $f(\omega)\leq f(\omega')$. By a standard monotone class argument, Inequality \eqref{FKG-inequality} is equivalent to
\begin{equation}
\E\big[f(\omega)g(\omega)\big] \geq \E\big[f(\omega)\big] \E\big[g(\omega)\big]
\end{equation}
for every two increasing $f,g\in L^2(\P)$.

\subsection{General theory of kernels.}\label{section-Kernels} Let $X$ be a set. Recall the following definitions:
\begin{itemize}
\item[$\bullet$] A map $k\colon X \times X\to\R$ is a {\bf positive definite kernel} if the matrix $[k(x,y)]_{x,y\in F}$ is positive for every finite subset $F\subset X$, i.e.\ 
\begin{equation} \label{def-PD}
\sum_{i,j=1}^n \overline{a_i}a_j k(x_i,x_j)\geq0
\end{equation}
for every $x_1,\ldots,x_n\in X$, $a_1,\ldots,a_n\in\C$. It is {\bf normalized} if $k(x,x)=1$ for every $x\in X$.

\vspace{1mm}

\item[$\bullet$] A map $k\colon X \times X\to[0,\infty)$ is a

\vspace{1mm}

\begin{itemize}
\item[$\diamond$] {\bf conditionally negative definite kernel} if $k(x,x)=0$ and $k(x,y)=k(y,x)$ for every $x,y\in V$ and
\begin{equation} \label{def-CND}
\sum_{i,j=1}^n a_ia_j k(x_i,x_j)\leq 0
\end{equation}
for every $x_1,\ldots,x_n\in V$ and $a_1,\ldots,a_n\in\R$ with $\sum_{i=1}^n a_i = 0$.

\vspace{1mm}

\item[$\diamond$] {\bf measure definite kernel}, if there exists a measure space $(\Omega,\skrib,\mu)$ and a map $S\colon X\to\skrib$, $x\mapsto S_x$, such that 
\begin{equation}
k(x,y)=\mu(S_x \Delta S_y)
\end{equation}
for every $x,y\in X$. 

\vspace{1mm}

\item[$\diamond$] {\bf $L^1$-kernel}, if there exists a measurable space $(Y,\nu)$ and a map $f\colon X\to L^1(Y,\nu)$ such that 
\begin{equation}
k(x,y) = \Vert f(x)-f(y) \Vert
\end{equation}
for every $x,y\in X$.
\end{itemize}
\end{itemize}

The reason for grouping the four concepts as above is that positive definiteness should be seen as opposite to the other three, which in turn are related as follows: clearly every measure definite kernel is an $L^1$-kernel. Indeed, take $f(x)=\1_{S_x}-\1_{S_{x_0}}\in L^1(\Omega,\mu)$ for a fixed base point $x_0\in X$. In fact, the two concepts agree by the following classical argument, see for instance \cite[Proposition 2.6]{CSV12}, which we include for the convenience of the reader.

\begin{prop}[Measure definite and $L^1$-kernels coincide] \label{prop-L1-MD}
Let $X$ be a countable set and let $k\colon X\times X\to[0,\infty)$ be an $L^1$-kernel. Then $k$ is measure definite.
\end{prop}

\proof 
The class of measure definite kernels is a convex cone, as can be verified directly using the disjoint union of two measure spaces. Indeed, let $k$ and $k'$ be measure definite kernels on $X$ and let $\lambda,\lambda'\ge0$. We claim that $k_0(x,y):=\lambda k(x,y)+\lambda'k(x,y)$
is measure definite. To see this, choose measure spaces $(\Omega,\mathcal B,\mu)$, $(\Omega',\mathcal B',\mu')$ together with maps $S:X\to\mathcal B$, $S':X\to\mathcal B'$ with the property that $k(x,y)=\mu(S_x\Delta S_y)$ and $k'(x,y)=\mu'(S_x'\Delta S_y')$. Without loss of generality, $\Omega$ and $\Omega'$ are disjoint sets (otherwise, consider $\Omega\times\{1\}$ and $\Omega'\times\{2\}$). We set $\Omega_0:=\Omega\sqcup\Omega'$ and $\mathcal B_0:=\{A\sqcup A':A\in\mathcal B, A'\in\mathcal B'\}$, which defines a $\sigma$-field on $\Omega_0$. For $A_0\in\mathcal B_0$, set
$$
\mu_0(A_0):=\lambda \mu(A_0 \cap \Omega) +\lambda' \mu'(A_0 \cap \Omega'),
$$
which yields a well-defined measure on $(\Omega_0,\mathcal B_0)$. Let $S^{(0)} \colon X \to \mathcal B_0 \, , \, S^{(0)}_x:=S_x \sqcup S_x'$. One checks that
$$
S^{(0)}_x \Delta S_y^{(0)} = (S_x\Delta S_y) \sqcup (S_x'\Delta S_y'),
$$
which implies that $k_0(x,y) = \mu_0( S_x^{(0)} \Delta S_y^{(0)})$, i.e.~measure definiteness of $k_0$. 

The class of measure definite kernels is also closed under pointwise convergence, see Proposition~\ref{prop-MD-PW-convergence} below.
The convex cone of $L^1$-kernels, on the other hand, is generated by {\bf cut-metrics}, i.e.~pull-backs of $\{0,1\}$-valued metric on two elements, see for instance \cite[Section 4.2]{DL97}. Since every cut-metric is measure-definite, the claim follows. \eproof

We single out the following standard fact about measure definite kernels used in the above proof due to its importance.

\begin{prop}[Pointwise limits of measure definite kernels, cf.~{\cite[Proposition 1.3]{RS98}}] \label{prop-MD-PW-convergence}
Let $X$ be a countable set. The class of measure definite kernels on $X$ is closed under pointwise convergence.
\end{prop}

A sketch of proof of the above proposition is provided below Lemma \ref{lemma-CSV}.

Every measure definite kernel is conditionally negative definite, but not conversely
, see Proposition~1.1 and the discussion thereafter in~\cite{RS98}.
Conditionally negative definite kernels have some particularly useful features: most importantly, $k\colon X\times X\to[0,\infty)$ is conditionally negative definite if and only if it is of the form 
\begin{equation} \label{equ-CND-HS}
k(x,y) = \Vert f(x)-f(y)\Vert_H^2
\end{equation} 
for a Hilbert space $H$ and a map $f\colon X\to H$, see e.g.~\cite[Theorem C.2.3]{BdlHV08}. This is equivalent, by a well-known result of Schoenberg, to
\begin{equation} \label{equ-Schoenberg}
k_\lambda \colon X\times X \to [0,1] \, , \, k_\lambda(x,y):=\e^{-\lambda k(x,y)}
\end{equation}
defining a positive definite kernel for every $\lambda\geq0$, see e.g.~\cite[Theorem C.3.2]{BdlHV08}.

Another important example of positive definite kernels, which explains the connection with percolation, is provided by the following lemma.

\begin{lemma}[Two-point functions are positive definite] \label{lemma-TPF-PD}
Let $\P$ be a general bond percolation on a graph $G=(V,E)$. Then the two-point function defines a normalized, positive definite kernel on $V\times V$.
\end{lemma}

\proof This observation is due to Aizenman and Newman \cite{AN84}. We include the short proof for the convenience of the reader: First note that $\tau$ is normalized as $\tau(v,v)=\P\big[v\leftrightarrow v\big]=1$ for every $v\in V$. Now let $v_1,\ldots,v_n\in V$ and $a_1,\ldots,a_n\in\C$. Let $\omega$ denote a random variable with law $\P$ and denote by $\skric$ the set of clusters in $\omega$. Then
\begin{flalign*}
\sum_{i,j=1}^n \overline{a_i}a_j\tau(v_i,v_j) &= \sum_{i,j} \overline{a_i}a_j \E\big[ \1_{\{v_i \leftrightarrow v_j\}}\big] = \E \biggl[  \ \sum_{i,j} \overline{a_i}a_j \1_{\{v_i \leftrightarrow v_j\}} \biggr] \\
& =  \E \biggl[ \ \sum_{C\in\skric} \, \sum_{\{v_i,v_j\} \subset C} \overline{a_i}a_j \biggr] = \E \biggl[ \ \sum_{C\in\skric} \, \Bigl| \sum_{v_i \in C} a_i \Bigr|^2 \biggr] \geq 0.
\end{flalign*}
Hence $\tau$ defines a positive definite kernel. \eproof

\section{Measure definite kernels, two-point functions and wall structures on general graphs}\label{mainsec MD2wall}

In this section, the key technical ingredient that $1-\P\big[u\leftrightarrow v]$ defines a measure definite kernel is proven in Lemma \ref{lemma-TPF-MD} in Section~\ref{sec MD2PWall}. Section~\ref{section-ConstructionWalls} contains the measure-theoretic construction regarding walls on general graphs needed in the sequel.

\subsection{Measure definite kernels and two-point functions.}\label{sec MD2PWall} In the setting of Lemma \ref{lemma-TPF-PD}, positive definiteness together with a standard result \cite[Proposition C.2.4(iii)]{BdlHV08} imply that
\begin{align*}
k\colon V\times V \to [0,1] \, , \, k(u,v) & =\P\big[u\not\leftrightarrow v\big] = 1-\P\big[u\leftrightarrow v \big] \\
&= \P\big[u\leftrightarrow u \big] - \P\big[u\leftrightarrow v \big]
\end{align*}
is a conditionally negative definite kernel. In fact, this definition actually provides examples of measure definite kernels as shown by the following lemma.

\begin{lemma}[Measure definite kernels and two-point functions]\label{lemma-TPF-MD}
Let $\P$ be a general bond percolation on a graph $G=(V,E)$. Then $\P\big[u\not\leftrightarrow v\big]$ defines a measure definite kernel on $V\times V$.
\end{lemma}

\proof 
To prove that the kernel
\begin{equation*}
k\colon V\times V \to [0,1] \, , \, k(u,v):=\P\big[u\not\leftrightarrow v\big]
\end{equation*}
is measure definite, it suffices to show that its restriction $k_{|F\times F}$ to every finite subset $F\subset V$ is measure definite (this follows e.g.~from Proposition~\ref{prop-MD-PW-convergence} or from~\cite[Corollary~6.17]{CDH10}).

Consider an arbitrary finite subset $F=\{v_1,\ldots,v_N\}\subset V$. Let $\omega$ denote the configuration of $\P$ and, for $v\in V$, let $C(v)$ denote the $\omega$-cluster of $v$ viewed as a subset of~$V$. Define random subsets $C_1,\ldots,C_N$ of $F$ as follows: set $C_1:= C(v_1)\cap F$ and, iteratively for $i\in\{2,\ldots,N\}$,
$$
\begin{aligned}
C_i:= 
\begin{cases}
C(v_i)\cap F \quad & \mbox{if}\quad v_i\notin\bigcup_{j=1}^{i-1}C_j,  \\
\emptyset\quad & \mbox{else.}
\end{cases}
\end{aligned}
$$
We obtain random, pairwise disjoint subsets $C_1,\ldots,C_N$ of $F$ which are either empty or are the restriction of one of the clusters of $\omega$ to the set $F$. Let $\mu_i$ denote the law of $C_i$ viewed as a probability measure on $\{0,1\}^F$. For $v\in F$, define
\begin{equation}
    S_v=\{\xi\in\{0,1\}^F \colon \xi(v)=1\}.
\end{equation} 
By definition we have that a pair of vertices $u,v\in F$ is not connected in $\omega$ if and only if there exists a unique $i\in\{1,\ldots,N\}$ such that $C_i$ contains $u$ and does not contain $v$, namely the $C_i$ which coincides with the restriction of $C(u)$ to $F$. It follows that for every $u,v\in F$,
\begin{flalign}
k(u,v) & = \P \Big[ \, u \in C_i \, \, \text{and} \, \, v \notin C_i \, \, \mbox{for a unique} \, \, i \in \{1,\ldots,N\} \, \Big] \nonumber \\
& = \sum_{i=1}^N \P\big[u \in C_i, v\notin C_i\big] \nonumber \\
& = \sum_{i=1}^N \mu_i(S_u \setminus S_v).  \label{eq-mdkernel1} 
\end{flalign}
By symmetry, we also have that 
\begin{equation}\label{eq-mdkernel2}
     k(u,v) = \sum_{i=1}^N \mu_i(S_v \setminus S_u)
\end{equation}  
for every $u,v\in F$. Summing up (\ref{eq-mdkernel1}) and (\ref{eq-mdkernel2}), we obtain that
\begin{equation} \label{eq-mdkernel3}
k(u,v)=\frac{1}{2} \sum_{i=1}^N \mu_i(S_u \Delta S_v) 
\end{equation}
for every $u,v\in F$. 
We see that $k_{|F\times F}$ is a convex combination of measure definite kernels and hence measure definite.
\eproof

\begin{remark} \label{rem-TPF-MD-background} 
In the context of group-invariant percolation on Cayley graphs, Lemma~\ref{lemma-TPF-MD} was obtained by the present authors earlier in a preliminary version of~\cite{MR23}. After completion, we realized that the result does not rely on the underlying group structure and, moreover, the general version can be leveraged to characterize the $L^1$-compression exponent, significantly improving similar results in the invariant situation (see Remark \ref{rem-L1andL1NE} for details). Lemma \ref{lemma-TPF-MD} thus fits into the general setting considered in this paper, and consequently has been moved here. \nopagebreak {\hfill\rule{2mm}{2mm}}
\end{remark}

\subsection{Wall structures on general graphs.}\label{section-ConstructionWalls} The ideas behind the measure theoretic construction described in this section go back to Robertson and Steger \cite{RS98} and were further developed and applied in \cite{CMV04,CDH10,CSV12}. For conciseness, we only recall the concepts strictly needed and refer to these works for details, background and related results.

Let $X$ be a countable set and denote by 
\begin{equation}
\Omega_X:=\{0,1\}^X\setminus\big\{(0,0,\ldots),(1,1,\ldots)\big\}
\end{equation}
the set of non-trivial subsets of $X$. We equip $\Omega_X$ with the $\sigma$-algebra $\skrib_X$ of Borel sets. For $x\in X$, we define 
\begin{equation}
S_x:= \big\{c\in\Omega_X\colon c(x)=1 \big\}\in\skrib_X
\end{equation}
to be the set of non-trivial subsets of $X$ containing $x$. We observe that $S_x\Delta S_y$ consists of all nontrivial subsets of $\Omega_X$ which {\bf separate} $x$ and $y$ in the sense that they contain exactly one of the two. Such sets may be interpreted as {\em walls} between $x$ and $y$ and this interpretation leads to the concept of {\em space with measured walls} \cite{CMV04}, see also \cite{CSV12}. For the purposes of this paper (i.e.~proving our two main results), it suffices to recall the following two technical lemmas. The first result shows that the space $(\Omega_X,\skrib_X)$ is a universal measurable space for measure definite kernels, resp.~$L^1$-kernels.

\begin{lemma}[Wall structures associated to measure definite and $L^1$-kernels] \label{lemma-CSV} Let $X$ be a countable set and let $k\colon X\times X\to[0,\infty)$ be a measure definite kernel or an $L^1$-kernel. Then there exists a regular Borel measure $\mu_k$ on $(\Omega_X,\skrib_X)$ such that 
\begin{equation}
k(x,y)=\mu_k(S_x\Delta S_y) 
\end{equation}
for every $x,y\in X$.
\end{lemma}

\noindent{\em Sketch of proof.}  We may assume that $k$ is measure definite by Proposition \ref{prop-L1-MD}. Thus, there exists a measure space $(\Omega,\skrib,\nu)$ and a map $T\colon X\to\skrib, x\mapsto T_x$ such that $k(x,y)=\nu(T_x\Delta T_y)$ for every $x,y\in X$. Without loss of generality, $\Omega$ does not contain any points belonging to all, resp.~none, of the sets $T_x$. Define a map 
$$
\phi\colon \Omega\to\Omega_X \, , \, \phi(\omega)(x):= \1_{T_x}(\omega).
$$
Since $\phi$ is measurable, the push-forward of $\nu$ under this map is a well-defined measure $\mu$ on $(\Omega,X,\skrib)$. Then
$$
\mu(T_x\Delta T_y)= \nu\big( \phi^{-1}(T_x\Delta T_y)\big) = \nu(T_x \Delta T_y)
$$
and we obtain a measure $\mu$ on $(\Omega_X,\skrib_X)$ implementing the measure definite kernel $k$. The fact that $\mu$ is regular may be concluded from the two properties that open subsets of $\Omega_X$ are $\sigma$-compact and that $\mu(K)<\infty$ for compact subsets $K\subset\Omega$. We refer to \cite[Proposition 1.2]{RS98} for details.
\eproof

\noindent{\em Sketch of proof of Proposition \ref{prop-MD-PW-convergence}.} For the convenience of the reader, we include a brief sketch of the proof and refer to \cite[Proposition 1.3]{RS98} for details. By Lemma \ref{lemma-CSV}, it suffices to prove that if $k$ is the pointwise limit of a sequence $(k_n)$ of measure definite kernels of the form $k_n(x,y)=\mu_n(S_x\Delta S_y)$, where each $\mu_n$ is a regular Borel measure on $(\Omega_X,\skrib_X)$, then $k$ is measure definite. By pointwise convergence, $\sup_n \mu_n(S_x\Delta S_y)<\infty$ for every $x,y\in X$. By compactness and the fact that $X$ is countable, we may assume that ${\mu_n}_{|S_x\Delta S_y}$ converges weakly for every $x,y\in X$. Using this observation, it is not difficult to conclude that $\Omega_X$ decomposes into a countable disjoint union of sets $T_i$ such that ${\mu_n}_{|T_i}$ converges weakly to $\mu^{(i)}$ for each $i\geq1$ and such that every $S_x\Delta S_y$ is contained in the union of at most finitely many of the sets $T_i$. Defining $\mu=\sum_{i=1}^\infty \mu^{(i)}$, the fact that $S_x\Delta S_y$ is contained in at most finitely many of the sets $T_i$ implies that
\begin{align*}
\mu(S_x \Delta S_y) & = \sum_{i=1}^\infty \mu^{(i)}\big(T_i \cap (S_x\Delta S_y)\big) = \lim_{n\to\infty} \sum_{i=1}^\infty \mu_n\big(T_i \cap (S_x\Delta S_y)\big) \\
& = \lim_{n\to\infty} \mu_n(S_x\Delta S_y) = k(x,y),
\end{align*}
i.e.~the kernel $k$ is indeed measure definite.\eproof

\begin{lemma}[Wall structures associated to negative definite kernels] \label{lemma-RS} Let $X$ be a countable set and let $k\colon X\times X\to[0,\infty)$ be a conditionally negative definite kernel. Then there exists a regular Borel measure $\mu_k$ on $(\Omega_X,\skrib_X)$ such that 
\begin{equation}
\sqrt{k(x,y)}=\mu_k(S_x\Delta S_y)
\end{equation}
for every $x,y\in X$.
\end{lemma}

\noindent{\em Sketch of proof.} We include a brief sketch of the proof and refer to \cite[Proposition 1.4]{RS98} for details. By taking limits along increasing unions of a sequence of finite sets (see \cite[Proposition 1.3 and 1.4]{RS98} for this part of the argument), it suffices to prove the lemma for finite $X$. In the finite case, recall from \eqref{equ-CND-HS} that we can represent
$$
k(x,y)=\Vert f(x)-f(y) \Vert_H^2
$$
for a map $f\colon X\to H$ into a, without loss of generality real, Hilbert space $H$. Replacing $H$ by the linear subspace generated by $\{f(x)\colon x\in X\}$ if necessary, we may assume that $H$ is finite dimensional. The homogeneous space $\Omega$ of all half-spaces of $H$ carries a natural measure $\mu$ which is invariant under (the unimodular group of) rigid motions of $H$ and can be normalized such that 
$$
\mu\big(\{ E\in\Omega \colon \xi\in E, \eta\notin E\} \big) = \Vert \xi-\eta \Vert_H
$$
for every $\xi,\eta\in H$. Defining $T_x:=\{E\in\Omega\colon f(x)\in E\}$, we obtain that
$$
\sqrt{k(x,y)}= \Vert f(x)-f(y) \Vert = \mu(S_x\setminus S_y) = \mu(S_y\setminus S_x),
$$
and thus $\sqrt{k}$ is indeed measure definite and thus of the desired form by Lemma \ref{lemma-CSV}.
\eproof

Note that Lemma \ref{lemma-RS} in particular implies that the square root of a conditionally negative definite kernel is measure definite. 

In the case that $X$ is the vertex set of a graph $G=(V,E)$, we will use the notation
\begin{equation}
W_e := S_x \Delta S_y \quad \mbox{if} \, \, e = [x,y]\in E,
\end{equation}
for the set of walls separating the endpoints of the edge.


\section{Coarse embeddability and $L^1$-compression exponent of graphs} \label{section-CEandL1}

The goal of this section is to prove Theorem \ref{theorem-CE} and Theorem \ref{theorem-L1NE} (and therefore Corollary \ref{cor-CE-groups} and Corollary~\ref{cor-L1NE-groups}), stated again below.

\begin{theorem}[Percolation and coarse embeddability into a Hilbert space for graphs] \label{theorem-Coarse-Embeddability}
Let $G=(V,E)$ be a locally finite, connected graph. Then $G$ admits a coarse embedding into a Hilbert space if and only if for every $p<1$, there exists a general bond percolation $\P$ with $\P\big[e\in\omega\big]>p$ for every $e\in E$ and such that the two-point function vanishes at infinity, i.e.
\begin{equation*} 
\lim_{r\to\infty} \sup \Big\{\P\big[u\leftrightarrow v\big] \colon u,v\in V, d(u,v)>r \Big\} = 0.
\end{equation*}
\end{theorem}

\medskip

\begin{theorem}[Percolation and the $L^1$-compression exponent for graphs] \label{theorem-L1-Compression-Exponent} 
Let $G=(V,E)$ be a locally finite, connected graph. Then $\alpha_1^*(G)$ is the supremum over those $\alpha\in[0,1]$ for which there exists $C>0$ such that for every $p<1$, there exists a general bond percolation $\P$ with $\P\big[e\in\omega\big]>p$ for every $e\in E$ and which satisfies the stretched exponential decay
\begin{equation*}
\e^{-\beta d(u,v)} \leq \P\big[u \leftrightarrow v \big] \leq \e^{-\gamma d(u,v)^\alpha} 
\end{equation*}
for every $u,v\in V$ and some $\beta,\gamma>0$ with $\beta/\gamma\leq C$.
\end{theorem}

In fact, we actually prove the following general quantitative result,

\begin{theorem}[Characterization of $L^1$-compression functions]\label{theorem-L1-Compression-Functions} Let $G=(V,E)$ be a locally finite, connected graph and let $\rho\colon[0,\infty]\to[0,\infty]$ be a function with $\{\rho=0\}=\{0\}$. Then the following are equivalent:
\begin{itemize}
\item[{\rm(i)}] There exists a Lipschitz function $f\colon V\to L^1$ satisfying
\begin{equation*}
\Vert f(u)-f(v) \Vert \geq c \rho( d(u,v))
\end{equation*}
for every $u,v\in V$ and a uniform constant $c>0$.
\item[{\rm(ii)}] There exists $C>0$ such that for every $p<1$, there exists a general bond percolation $\P$ with $\P\big[e\in\omega\big]>p$ for every $e\in E$ and with
\begin{equation*}
\e^{-\beta d(u,v)} \leq \P\big[u \leftrightarrow v \big] \leq \e^{-\gamma \rho(d(u,v))} 
\end{equation*}
for every $u,v\in V$ and some $\beta,\gamma>0$ with $\beta/\gamma\leq C$.
\end{itemize}
\end{theorem}


\medskip

\begin{remark}[Lipschitz constants and the two-point function decay rate] In fact, in Theorem \ref{theorem-L1-Compression-Functions}, we show that if $f\colon V\to L^1$ is $L$-Lipschitz with $\Vert f(u)-f(v) \Vert \geq c \rho( d(u,v))$, then $\beta/\gamma\leq  L/c$ in (ii). Conversely, if $\beta/\gamma\leq C$, then there exists an $C$-Lipschitz function $f\colon V\to L^1$ with $\Vert f(u)-f(v) \Vert \geq \rho( d(u,v))$, see the proof of Theorem \ref{theorem-L1-Compression-Functions} for details.\nopagebreak {\hfill\rule{2mm}{2mm}}
\end{remark}

We also show that the percolation condition from  \cite{MR23} characterizing property (T) (see~Section~\ref{section-ClosingRemarks}~(A) for the precise statement of the condition) trivializes for general percolations in the sense that no locally finite, infinite, connected graph exists which satisfies that condition.

\begin{prop}[Connectivity decay and positive associations on infinite graphs] \label{prop-GeneralDecay}
Let $G=(V,E)$ be a locally finite, connected, infinite graph. Then for every $p<1$, there exists a general bond percolation~$\P$ with $\P\big[e\in\omega\big]>p$ for every $e\in E$, which satisfies the FKG-inequality and has connectivity decay, i.e.\
\begin{equation}
\inf_{u,v\in V} \, \P\big[u \leftrightarrow v \big] = 0.
\end{equation}
\end{prop}


\subsection{Construction of one-parameter family of percolations.} \label{section-ConstructionPerc} To every wall structure on a graph, we associate in this section a $1$-parameter family of percolations with parameter $t>0$. Here the parameter is the intensity of a Poisson process on the wall structure and intuitively corresponds to large marginals when $t$ is small and small marginals when $t$ is large.

\begin{prop}[Percolations from negative definite kernels] \label{prop-GeneralConstruction} Let $G=(V,E)$ be a locally finite, connected graph and let $k\colon V \times V\to[0,\infty)$ be a conditionally negative definite kernel. Then for every $t>0$, there exists a general bond percolation $\P_t$ which satisfies the FKG-inequality and the following:
\begin{enumerate}
\item[{\rm(i)}] for every $[u,v]=e\in E$, 
\begin{equation} \label{equ-GeneralConstructionMarginals}
\P_t\big[e\in E\big] = \exp\big(-t \sqrt{k(u,v)} \big).
\end{equation}
\item[{\rm(ii)}] for every $u,v\in V$,
\begin{equation}  \label{equ-GeneralConstructionTPF}
\exp\big(- t\delta d(u,v) \big) \leq \P_t\big[u \leftrightarrow v\big] \leq \exp\big(-t \sqrt{k(u,v)}\big),
\end{equation}
where $\delta := \sup_{u\sim v} \sqrt{k(u,v)}$.
\end{enumerate}
Moreover, there exists a monotonically decreasing coupling of the family $(\P_t)_{t>0}$. 
\end{prop}

\proof 
Let $(\Omega_V,\skrib_V,\mu_k)$ be as guaranteed by Lemma~\ref{lemma-RS}.
Let $\eta$ be a Poisson process on $\Omega_V\times[0,\infty)$ with intensity measure $\mu_k\otimes\d s$, where $\d s$ denotes Lebesgue measure -- i.e.~$\eta$ is a random counting measure on $\Omega_V\times[0,\infty)$ such that
\begin{enumerate}
\item for every measurable $B\subset\Omega_V\times[0,\infty)$,
\begin{equation*}
\P\big[\eta(B)=i\big] = \e^{-(\mu_k\otimes\d s)(B)} \frac{(\mu_k\otimes\d s)(B)^i}{i!}
\end{equation*}
for every $i\in\{0,1,\ldots\}$.

\vspace{1mm}

\item for every $n\in\N$ and every collection of pairwise disjoint, measurable $B_1,\ldots,B_n\subset\Omega_V\times[0,\infty)$, the random variables $\eta(B_1),\ldots,\eta(B_n)$ are independent.
\end{enumerate}
For $t\in[0,\infty)$, define the configuration $\omega_t$ of a bond percolation on $G$ as follows: for each edge $[u,v]=e$, let $e\in\omega_t$ if and only if 
\begin{equation*}
\eta\big( W_e \times [0,t) \big)=0.
\end{equation*} 
Let $\P_t$ denote the law of $\omega_t$. Observing that $\omega_t\subset\omega_s$ for $s<t$, we see that the family $(\P_t)_{t>0}$ has the desired monotone coupling. The FKG-inequality for the Poisson process, see \cite[Theorem 20.4]{LP17}, implies that $\P_t$ satisfies the FKG-inequality (note that every increasing event for $\omega_t$ is a decreasing event for $\eta$). Moreover, 
\begin{equation*}
\P_t\big[e\in E\big] = \PP\Big[ \eta\big( W_e \times [0,t) \big)=0 \Big]  = \exp\big(-t \sqrt{k(u,v)} \big),
\end{equation*}
which shows (i). For the lower bound, let $u,v\in V$ and let $e_1,\ldots,e_{d(u,v)}$ be a path of edges joining $u$ to $v$. Note that 
\begin{equation*}
\Big\{ \eta\Big( \cup_{i=1}^{d(u,v)} W_{e_i} \times [0,t) \Big) = 0 \Big\} \subset \big\{ u\overset{\omega_t}{\longleftrightarrow} v\big\}
\end{equation*}
and thus 
\begin{equation*}
\P_t\big[u \leftrightarrow v\big] \geq \exp\Big( -t \mu_k\big(\cup_{i=1}^{d(u,v)} W_{e_i}\big)\Big) \geq \exp\big(- t\delta d(u,v) \big),
\end{equation*}
where $\delta := \sup_{u\sim v} \sqrt{k(u,v)}$. Finally, let $u,v\in V$ and suppose that 
\begin{equation*}
\eta\big((S_u\Delta S_v)\times [0,t)\big)>0.
\end{equation*}
Then there exists $(c,s)\in\supp(\eta)$ which satisfies $c(u)\not=c(v)$ and $s<t$. Now consider any path joining $u$ to $v$ in $G$ and label the appearing edges $e_1,e_2,\ldots,e_n$. Write $e_i=[u_i,v_i]$ and $u=u_0$. Then we may choose $i\in\{0,\ldots,n-1\}$ minimal such that $c(e_i)\neq c(e_{i+1})$. In particular, $(c,s)\in (S_{u_i}\Delta S_{u_{i+1}})\times [0,t)$ and hence 
\begin{equation*}
\eta\big( W_{e_i} \times [0,t)\big)>0.
\end{equation*} 
It follows that this path is not contained in $\omega_t$. Since the path was arbitrary, we obtain that
\begin{equation*}
\Big\{ \eta\big((S_u\Delta S_v)\times [0,t)\big)>0 \Big\} \subset \big\{ u\overset{\omega_t}{\longleftrightarrow} v\big\}^c,
\end{equation*}
and hence
\begin{equation*}
\P_t \big[u \leftrightarrow v] \leq \PP \Big[\eta\big((S_u\Delta S_v)\times [0,t)\big)=0\Big] = \exp\big(-t \sqrt{k(u,v)}\big),
\end{equation*} 
which proves (ii) and thus completes the proof.\eproof

The same arguments as in the above proof of Proposition \ref{prop-GeneralConstruction} yield the following.

\begin{cor}[Percolations from measure definite and $L^1$-kernels] \label{cor-ConstructionL1} Let $G=(V,E)$ be a locally finite, connected graph and let $k\colon V \times V\to[0,\infty)$ be a measure definite kernel or an $L^1$-kernel. Then for every $t>0$, there exists a general bond percolation $\P_t$ which satisfies the FKG-inequality and the following:
\begin{enumerate}
\item[{\rm(i)}] for every $[u,v]=e\in E$, 
\begin{equation} \label{equ-L1ConstructionMarginals}
\P_t\big[e\in \omega \big] = \exp\big(-t k(u,v) \big).
\end{equation}
\item[{\rm(ii)}] for every $u,v\in V$, 
\begin{equation} \label{equ-L1ConstructionTPF}
\exp\big(- t\delta d(u,v) \big) \leq \P_t\big[u \leftrightarrow v\big] \leq \exp\big(-t k(u,v) \big),
\end{equation}
where $\delta := \sup_{u\sim v} k(u,v)$.
\end{enumerate}
Moreover, there exists a monotonically decreasing coupling of the family $(\P_t)_{t>0}$. 
\end{cor}

\subsection{Proofs of the main results.} \label{section-proofs} With the percolation construction established, we now proceed to the proofs of the main results (Theorem \ref{theorem-Coarse-Embeddability}- Theorem \ref{theorem-L1-Compression-Functions} and Proposition \ref{prop-GeneralDecay}).

\medskip
\noindent{\bf Proof of Theorem \ref{theorem-Coarse-Embeddability}:} ``$\Leftarrow:$" Suppose that for every $p<1$, there exists a general bond percolation $\P$ with $\P\big[e\in\omega\big]>p$ for every $e\in E$ and such that the two-point function tends to zero off tubes. Choose a sequence $(\P_n)$ of general bond percolations such that $\P_n\big[e\in\omega\big]>1-1/n$ for every $e\in E$ and such that the corresponding two-point function 
\begin{equation*}
\tau_n \colon V\times V\to[0,1]\, , \, \tau_n(u,v)=\P_n\big[u\leftrightarrow v\big]
\end{equation*}
tends to zero off tubes. By Lemma \ref{lemma-TPF-PD}, $\tau_n$ is a normalized, positive definite kernel. Let $R,\eps>0$. We now verify that $\tau_n$ has $(R,\eps)$-variation for $n$ large enough: if $\gamma$ is any path in $G$ of length at most $R$, then it follows that
\begin{equation*}
\P_n\big[\gamma\subset\omega \big] \geq 1-R\Big(1-\inf_{e\in E} \P_n\big[e\in\omega\big] \Big) \geq 1-R/n
\end{equation*}
and thus
\begin{equation}\label{statement coarse}
\lim_{n\to\infty} \inf \big\{ \tau_n(u,v) \colon d(u,v)\leq R \big\}=1.
\end{equation}
To show that $G$ admits a coarse embedding into a Hilbert space, we make use of the implication ``(iii) $\Rightarrow$ (i)" of the following general fact: for a locally finite, connected graph $G=(V,E)$ the following are equivalent:
\begin{enumerate}
\item[{\rm (i)}] $G$ admits a coarse embedding into a Hilbert space;

\vspace{1mm} 

\item[{\rm (ii)}] there exists a conditionally negative definite kernel $k\colon V\times V\to [0,\infty)$ together with maps $\rho_1,\rho_2\colon[0,\infty)\to[0,\infty)$ such that $\rho_1(t)\to\infty$ as $t\to\infty$ and 
\begin{equation}
\rho_1(d(u,v))\leq k(u,v) \leq \rho_2(d(u,v))
\end{equation}
for every $u,v\in V$;

\vspace{1mm} 

\item[{\rm (iii)}] for every $R,\eps>0$, there exists a normalized, positive definite kernel $k\colon V\times V\to\R$ which tends to zero off tubes, i.e.
\begin{equation}
\lim_{r\to\infty} \sup \big\{k(u,v) \colon d(u,v)\geq r \big\} = 0,
\end{equation}
and has $(R,\eps)$-variation, i.e.~$d(u,v)\leq R$ implies $|1-k(u,v)|<\eps$.
\end{enumerate}

Indeed, the equivalence between (i) and (ii) can easily be deduced from the fact, mentioned already in \eqref{equ-CND-HS}, that $k\colon V\times V\to[0,\infty)$ is a conditionally negative definite kernel if and only if $k(u,v)=\Vert f(u)-f(v)\Vert_H^2$ for a Hilbert space $H$ and map $f\colon V\to H$. 

The implication from (ii) to (iii) follows from Schoenberg's theorem, stated in \eqref{equ-Schoenberg}, which asserts that if $k\colon V\times V\to[0,\infty)$ is a conditionally negative definite kernel, then $e^{-\lambda k}$ is a positive definite kernel for every $\lambda\geq 0$.

Finally, the fact that (iii) implies (ii) may be seen by the following construction: suppose (iii) holds and find for every $n\geq1$ some $R_n\geq \max\{n,R_{n-1}\}$ and a normalized, positive definite kernel $k_n$ such that $d(u,v)\leq n$ implies $|k_n(u,v)-1|<2^{-n}$ and such that $d(u,v)>R_n$ implies $k_n(u,v)<1/2$. It can then be checked that $k(u,v):=\sum_{n=1}^\infty 1-k_n(u,v)$ defines a conditionally negative definite kernel with $k(u,v)\leq2d(u,v)+1$ and $k(u,v)\geq 2^{-1}Q(d(u,v))$, where $Q(t):= \min\{n\geq1 \colon t \leq R_n\}$. We refer to \cite[Theorem 3.2.8]{W09} or \cite[Theorem 11.16]{R03} for details.

Thus, \eqref{statement coarse} and the implication ``(iii) $\Rightarrow$ (i)" of the above fact imply that $G$ admits a coarse embedding into a Hilbert space.

\medskip

\noindent$``\Rightarrow:"$ We now prove the converse direction and assume that $G$ admits a coarse embedding into a Hilbert space. By  the implication ``(i)$\Rightarrow$(ii)" of the aforementioned fact, there exists a conditionally negative definite kernel $k\colon V\times V\to [0,\infty)$ and maps $\rho_1,\rho_2\colon[0,\infty)\to[0,\infty)$ such that $\rho_1(t)\to\infty$ as $t\to\infty$ and 
\begin{equation} \label{equ-proofCEkBound}
\rho_1(d(u,v))\leq k(u,v) \leq \rho_2(d(u,v))
\end{equation}
for every $u,v\in V$. By Proposition \ref{prop-GeneralConstruction}, there exists for every $t>0$ a general bond percolation $\P_t$ with 
\begin{equation} \label{eq-proofCEEdgeBound}
\P_t\big[e\in E\big] = \exp\big(-t \sqrt{k(u,v)} \big)
\end{equation}
for every $[u,v]=e\in E$ and with
\begin{equation} \label{eq-proofCETPFBound}
\P_t\big[u \leftrightarrow v\big] \leq \exp\big(-t \sqrt{k(u,v)}\big)
\end{equation}
for every $u,v\in V$. Using \eqref{eq-proofCEEdgeBound} and the upper bound in \eqref{equ-proofCEkBound}, we see that
\begin{equation} \label{eq-CoarseEmbed1}
\inf_{e\in E} \P_t\big[e\in E\big]  = \exp\Big(-t \sup_{e=[u,v]\in E} \sqrt{k(u,v)} \Big) \geq \e^{-t\sqrt{\rho_2(1)}} \longrightarrow 1 \quad \mbox{as} \ t \to 0,
\end{equation}
and it follows that $\P_t$ has marginals at least $p$ for $t$ sufficiently small. Moreover, for any $t>0$, \eqref{eq-proofCETPFBound} and the lower bound in \eqref{equ-proofCEkBound} imply that
\begin{equation}
\limsup_{r\to\infty} \ \sup \Big\{\P\big[u\leftrightarrow v\big] \colon u,v\in V, d(u,v)>r \Big\} \leq \limsup_{r\to\infty} \, \e^{-t\sqrt{\rho_1(r)}}=0,
\end{equation}
i.e.~the two-point function of $\P_t$ vanishes at infinity. \eproof

\medskip
\noindent{\bf Proof of Theorem \ref{theorem-L1-Compression-Exponent}:} Let $f\colon V\to L^1$ be an $L$-Lipschitz function such that 
\begin{equation} \label{equ-L1Compr1}
\Vert f(u)-f(v) \Vert \geq c d(u,v)^\alpha
\end{equation}
for every $u,v\in V$ and some uniform constant $c>0$. Corollary \ref{cor-ConstructionL1} applied with the $L^1$-kernel
\begin{equation*}
k \colon V\times V \to [0,\infty) \, , \, k(u,v) := \Vert f(u)-f(v) \Vert
\end{equation*} 
yields the existence of a family $(\P_t)_{t>0}$ of general bond percolations with 
\begin{equation*}
\P_t\big[e\in E\big] = \exp\big(-t k(u,v) \big)
\end{equation*}
for every $[u,v]=e\in E$ and
\begin{equation*}
\exp\big(- t\delta d(u,v) \big) \leq \P_t\big[u \leftrightarrow v\big] \leq \exp\big(-t k(u,v) \big) \leq  \exp\big( -tcd(u,v)^\alpha \big),
\end{equation*}
for every $u,v\in V$, where 
\begin{equation*}
\delta := \sup_{u\sim v} k(u,v) \leq L
\end{equation*} 
because $f$ is $L$-Lipschitz. Setting 
\begin{equation*}
\beta_t := tL \qquad \mbox{and} \qquad \gamma_t := tc, 
\end{equation*}
it immediately follows that the two-point function of $\P_t$ has stretched exponential decay in the sense that
\begin{equation*}
\e^{-\beta_t d(u,v)} \leq \P_t\big[u \leftrightarrow v \big] \leq \e^{-\gamma_t d(u,v)^\alpha}
\end{equation*}
for every $t>0$ with 
\begin{equation*}
\sup_{t>0} \, \frac{\beta_t}{\gamma_t} = \frac{L}{c} =: C < \infty.
\end{equation*} 
Moreover, the lower bound on the distortion in \eqref{equ-L1Compr1}, together with the same argument used in \eqref{eq-CoarseEmbed1} above, implies that $\P_t$ has marginals at least $p$ for $t$ sufficiently small.

Conversely, suppose there exists a constant $C>0$ and a sequence $(\P_n)$ of general bond percolations, such that $\P_n$ satisfies $\P_n\big[e\in\omega\big]>1-1/n$ for every $e\in E$ and the stretched exponential two-point function decay
\begin{equation} \label{equ-ProofL1Bounds}
\e^{-\beta_n d(u,v)} \leq \P_n\big[u \leftrightarrow v \big] \leq \e^{-\gamma_n d(u,v)^\alpha} 
\end{equation}
for every $u,v\in V$ and some $\beta_n,\gamma_n>0$ with $\beta_n/\gamma_n\leq C$. Note that since the marginals tend to $1$ as $n\to\infty$, we must have $\gamma_n\to0$ and thus also $\beta_n\to0$ as $n\to\infty$. Set
\begin{equation}
k_n \colon V\times V \to[0,\infty) \, , \, k_n(u,v):= \frac{\P_n\big[u\not\leftrightarrow v\big]}{\gamma_n}.
\end{equation}
By Lemma \ref{lemma-TPF-MD} and the obvious fact that scalar multiplication preserves measure definiteness, each $k_n$ defines a measure definite kernel. Moreover, the lower bound in \eqref{equ-ProofL1Bounds} and the fact that $\beta_n\to0$ imply that
\begin{align}
   \nonumber\limsup_{n\to\infty} \ k_n(u,v) & = \limsup_{n\to\infty} \ \frac{1-\P_n\big[u \leftrightarrow v \big]}{\gamma_n} \\
& \nonumber\leq \limsup_{n\to\infty} \ \frac{1-\e^{-\beta_n d(u,v)}}{\gamma_n}  \\
&= \nonumber\limsup_{n\to\infty} \ \frac{\beta_n}{\gamma_n} \ \frac{1-\e^{-\beta_n d(u,v)}}{\beta_n} \\
&\leq C d(u,v) \label{eq-kLip} 
\end{align}
for every $u,v\in V$. By countability, we may go over to a subsequence such that 
\begin{equation}
k \colon V\times V \to[0,\infty) \, , \, k(u,v):= \lim_{n\to\infty} k_n(u,v)
\end{equation}
exists. As a pointwise limit of measure definite kernels, $k$ is measure definite, see Proposition \ref{prop-MD-PW-convergence}. By \eqref{eq-kLip}, it is Lipschitz with constant $C$. On the other hand, the upper bound in \eqref{equ-ProofL1Bounds} and the fact that $\gamma_n\to0$ imply that
\begin{equation} \label{eq-kalphagrowth}
    k(u,v) = \lim_{n\to\infty}  k_n(u,v) \geq \liminf_{n\to\infty} \ \frac{1-\e^{-\gamma_n d(u,v)^\alpha}}{\gamma_n} = d(u,v)^\alpha
\end{equation}
for every $u,v\in V$. Finally, Proposition \ref{prop-L1-MD} implies that $k$ is an $L^1$-kernel. Combining this fact with \eqref{eq-kLip} and \eqref{eq-kalphagrowth} shows that $\alpha_1^*(G)\geq\alpha$. \eproof

\begin{remark}(Difference to the equivariant case) \label{rem-L1andL1NE} Let us re-emphasize the differences to the equivariant case below. The construction in Proposition \ref{prop-GeneralConstruction} differs from the construction in the proof of \cite[Theorem 3.6]{MR23} in applying to general conditionally negative definite kernels on graphs and in directly providing the monotone coupling using 
the auxiliary Poisson process on $\Omega_V\times[0,\infty)$.
The difference between the equivariant and the non-equivariant compression exponents manifests itself in the proof of the existence of the appropriate measure definite kernel compared to the proof  of \cite[{Theorem~5.8~``(ii)$\Rightarrow$(iii)''}]{MR23}. While Lemma~\ref{lemma-TPF-MD} can be used in both situations to produce from the percolations a measure definite kernel with growth at least $d^\alpha$, this is sufficient to conclude $\alpha_1^*\geq\alpha$ only in the non-equivariant case. In the equivariant case, although the constructed kernel inherits invariance from the invariant percolations, it is not known that there actually is an associated isometric action implementing the kernel (in the sense of \eqref{def-L1}, see Section \ref{section-ClosingRemarks} (B) for details) -- this is the same problem as the ``missing implication'', which was raised in \cite[Proposition 2.8 \&~Remark~2.9]{CSV12}. In general, this property is only known for kernels which are the {\em square-root} of a conditionally negative definite kernel, where it is due to Robertson and Steger \cite[Proof of Theorem 2.1]{RS98} (see also \cite[Theorem 6.25]{CDH10} for a generalization). Due to having to take the square root in the equivariant situation, the current argument only allows to conclude $\alpha_1^{\scriptscriptstyle{\#}} \geq\alpha/2$ there. However, note that the kernel, whose square-root is taken, is already measured definite -- this is stronger than the conditional negative definiteness that is required for this argument. \nopagebreak {\hfill\rule{2mm}{2mm}}
\end{remark}

\begin{remark}[Percolation and $L^p$-compression exponents] \label{rem-Lp}
Let $G=(V,E)$ be a locally finite, connected, infinite graph and let $p>1$. Similarly to \eqref{equ-DefL1Compression}, the {\bf $L^p$-compression exponent} $\alpha_p^*(G)$ of $G$ is defined as the supremum over those $\alpha\in[0,1]$ such that there exist a measure space $(Y,\nu)$ and a Lipschitz function $f\colon V\to L^p(\nu)$ satisfying 
\begin{equation}
\Vert f(u)-f(v) \Vert \geq c d(u,v)^\alpha
\end{equation}
for every $u,v\in V$ and a uniform constant $c>0$, cf.~\cite{GK04,NP11}. As a by-product of the proof of Theorem~\ref{theorem-L1-Compression-Functions}, we obtain that if there exists $C>0$ such that for every $q<1$, there exists a general bond percolation $\P$ with $\P\big[e\in\omega\big]>q$ for every $e\in E$ and with
\begin{equation} \label{equ-LpTPFDecay}
\e^{-\beta d(u,v)^p} \leq \P\big[u \leftrightarrow v \big] \leq \e^{-\gamma \rho^p(d(u,v))} 
\end{equation}
for every $u,v\in V$ and some $\beta,\gamma>0$ with $\beta/\gamma\leq C$, then there exists a measure definite kernel $k\colon V\times V\to[0,\infty)$ such that 
\begin{equation*}
\rho^p(d(u,v))\leq k(u,v) \leq Cd(u,v)^p
\end{equation*} 
for every $u,v\in V$ and a uniform constant $C>0$. It then follows e.g.~from \cite[Proposition 2.6]{CSV12} that there exists a Lipschitz function $f\colon V\to L^p$ satisfying
\begin{equation} \label{equ-LpBound}
\Vert f(u)-f(v) \Vert \geq c \rho( d(u,v))
\end{equation}
for every $u,v\in V$ and a uniform constant $c>0$ (in fact, one may choose $\Vert f(u)-f(v) \Vert_p =k^{1/p}(u,v)$). In particular, in the special case that $\rho(t)=t^\alpha$, the existence of general bond percolations with large marginals and two-point function decay as in \eqref{equ-LpTPFDecay} implies $\alpha_p^*(G)\geq\alpha$.
Regarding the converse direction, i.e.~starting with some Lipschitz function $f\colon V\to L^p$ satisfying \eqref{equ-LpBound}, the kernel $k(u,v):=\Vert f(u)-f(v)\Vert_p^p$ need not be measure definite (as seen for $p=2$ in Section \ref{section-Kernels}). On the other hand, for $p\in(1,2]$, it is known that $k$ is conditionally negative definite, see e.g.~\cite[Proposition 6.3]{CDH10}. In particular, Proposition \ref{prop-GeneralConstruction} yields the existence of percolations with arbitrarily large marginals and two-point function decay as in \eqref{equ-LpTPFDecay}, but with the function $\rho$ replaced by $\rho^{1/2}$ (since $p>1$). Thus, if we start with some Lipschitz function $f\colon V\to L^p$ satisfying \eqref{equ-LpBound} with $\rho(t)=t^\alpha$, then this leads only to exponent $\alpha/2$ in \eqref{equ-LpTPFDecay}. \nopagebreak {\hfill\rule{2mm}{2mm}}
\end{remark}

\noindent{\bf Proof of Theorem \ref{theorem-L1-Compression-Functions}:} This follows exactly as in the above proof of Theorem \ref{theorem-L1-Compression-Exponent} with a general function $\rho$ replacing the function $\rho(t)=t^\alpha$.\eproof

\noindent{\bf Proof of Proposition \ref{prop-GeneralDecay}:} Let $G=(V,E)$ be a locally finite, connected, infinite graph. For $n\in\N$, define $A_n:=\{n,n+1,\ldots\} \subset\N$ and observe that $|A_n\Delta A_{m}|=\max\{n,m\}-\min\{n,m\}$. Fix a root $o\in V$ and define for each vertex $v$, $|v| := d(o,v)$ and $T_v:= A_{|v|}.$ It follows that
\begin{equation}
k \colon V\times V \to [0,\infty)\, , \, k(u,v):= |T_u\Delta T_v|
\end{equation}
is a symmetric, normalized and unbounded measure definite kernel with
\begin{equation} \label{equ-GeneralDecay1}
\sup\big\{k(u,v) \colon d(u,v)\leq R\big\}\leq R.
\end{equation}
Let $(\P_t)_{t>0}$ be a family of percolations as shown to exist in Corollary~\ref{cor-ConstructionL1}. For every $p<1$, \eqref{equ-GeneralDecay1} and \eqref{equ-L1ConstructionMarginals} guarantee that $\P_t$ has marginals at least $p$ for $t$ sufficiently close to $0$. Moreover, $\P_t$ satisfies the FKG-inequality. Finally, \eqref{equ-L1ConstructionTPF} and unboundedness of $k$ imply that
$$
\inf_{u,v\in V} \P_t\big[u \leftrightarrow v \big]  \leq \exp \Big( - t \sup_{u,v \in V} k(u,v) \Big) = 0,
$$
i.e.~connectivity decay.
\eproof

\begin{remark}(Connectivity decay vs.~vanishing connectivity) The kernel $k$ in the above proof is unbounded, but not proper as $k(u,v)=0$ whenever $d(o,u)=d(o,v)$. Thus the two-point function constructed in the above proof does not necessarily vanish at infinity (of course $G$ also need not coarsely embed into a Hilbert space). In a similar direction, let us mention that one could easily define proper measure definite kernels on $G$, even with linear growth, e.g.~by considering the induced distance of a spanning tree (every tree distance is a measure definite kernel, see e.g.~\cite[p.~251, Remark~(i)]{RS98}). The problem with such a kernel, however, is that it does not necessarily satisfy a uniform upper bound of the form $\sup \big\{ k(u,v) \colon d(u,v) \leq r \big\}<\infty$
and hence would not yield large marginals for small intensities.\nopagebreak {\hfill\rule{2mm}{2mm}}
\end{remark}

\section{Closing remarks} \label{section-ClosingRemarks}

\begin{itemize}
\item[(A)] The two results described below, obtained recently by the present authors \cite{MR23}, showed that geometric properties beyond amenability can be characterized using percolation. These results were also the main inspirations for the present work. Recall that given a finitely generated group $\Gamma$ and a finite, symmetric generating set $S$, the {\bf Cayley graph} of $\Gamma$ with respect to $S$ is the graph $G=(V,E)$ with $V=\Gamma$ and edges $\{[g,gs]\colon g\in\Gamma,s\in S\}$. A~{\bf$\Gamma$-invariant bond percolation} is a general bond percolation, which is invariant under left translations. The group $\Gamma$ has the {\bf Haagerup property}, if there exists an affine isometric action of $\Gamma$ on a Hilbert space $H$ which is metrically proper (i.e.,~$\big|\{g\in\Gamma\colon gB\cap B\neq \emptyset\}\big|<\infty$ for every bounded $B\subset H$). We may also point out that this implies, by considering the orbit of a point, the existence of an equivariant map $f:\Gamma\to H$ such that $\inf \{\Vert f(g)-f(h)\Vert : d_G(g,h)\ge r\}\to \infty$ as $r\to\infty$, which thus defines an equivariant coarse embedding. The converse is also true, as follows e.g.~by noting that $\Vert f(g)-f(h)\Vert^2$ defines a proper, $\Gamma$-invariant, conditionally negative definite kernel and \cite[Theorem 2.1.1]{CCJJV01}. The group has {\bf property (T)} if every affine isometric action of~$\Gamma$ on a Hilbert space has a fixed point. 
An important fact is that $\Gamma$ has both properties if and only if it is finite. See \cite{BdlHV08,BO08,CCJJV01} for more information on these properties. 
We are now in a position to recall the interplay with invariant percolation\footnote{We point out that large marginals were equivalently expressed in \cite{MR23} by requiring large expected degree, i.e.~$\E\big[\deg_\omega(v)\big]>\alpha \, \deg(v)$ for every $v\in V$ and some $\alpha<1$ close to $1$.}:

\vspace{1mm}

\begin{itemize}
\item[$\bullet$] The finitely generated group $\Gamma$ has the Haagerup property if and only if some, equivalently every, Cayley graph $G=(V,E)$ has the property that for every $p<1$, there exists a $\Gamma$-invariant bond percolation $\P$ with  $\P\big[e\in\omega\big] \geq p$ for every $e\in E$ and such that the two-point function vanishes at infinity.

\vspace{1mm}

\item[$\bullet$] The finitely generated group $\Gamma$ has property (T) if and only if some, equivalently every, Cayley graph $G=(V,E)$ has the property that there exists $p^*<1$ such that every $\Gamma$-invariant bond percolation $\P$ with $\P\big[e\in\omega\big]>p^*$ for every $e\in E$ exhibits long-range order, i.e.~satisfies~$\inf_{u,v\in V} \P\big[u\leftrightarrow v\big] > 0$.
\end{itemize}

\vspace{1mm}

We point out that the direct implication in the second result above is due to \cite{LS99}. 

\vspace{1mm}

\item[(B)] For completeness, we now recall explicitly what is known in the equivariant setting. The {\bf equivariant $L^1$-compression exponent} $\alpha_1^{\scriptscriptstyle{\#}}(\Gamma)$ of a finitely generated group $\Gamma$ is the supremum over those $\alpha\in[0,1]$, such that there exists a measure space $(Y,\nu)$ and an affine isometric action of $\Gamma$ on $L^1(Y,\nu)$ together with a $\Gamma$-equivariant map $f\colon\Gamma\to L^1(Y,\nu)$ satisfying 
\begin{equation} \label{def-L1}
\Vert f(g)-f(h) \Vert \geq cd(g,h)^\alpha
\end{equation}
for every $g,h\in\Gamma$ and a uniform constanct $c>0$, where $d$ is some word  metric on $\Gamma$, cf.~\cite{GK04,NP11}. Towards a probabilistic understanding of this exponent,  it was shown in \cite{MR23} that for $\alpha\in[0,1]$, each of the following conditions implies the next (and ''(i)$\Leftrightarrow$(ii)'' for amenable $\Gamma$):

\vspace{1mm}

\begin{enumerate}
\item[{\rm(i)}] The equivariant $L^1$-compression exponent $\alpha_1^{\scriptscriptstyle{\#}}(\Gamma)$ is at least $\alpha$.

\item[{\rm(ii)}] Some, equivalently every, Cayley graph $G=(V,E)$ has the property that there exists $C>0$ such that for every $p<1$, there exists a $\Gamma$-invariant bond percolation $\P$ with $\P\big[e\in\omega\big]>p$ for every $e\in E$ and which satisfies the stretched exponential decay 
\begin{equation*}
\e^{-\beta d(u,v)} \leq \P\big[u \leftrightarrow v \big] \leq \e^{-\gamma d(u,v)^\alpha} 
\end{equation*}
for all $u,v\in V$ and some $\beta,\gamma>0$ with $\beta/\gamma\leq C$.

\item[{\rm(iii)}] The equivariant $L^1$-compression exponent $\alpha_1^{\scriptscriptstyle{\#}}(\Gamma)$ is at least $\alpha/2$.
\end{enumerate}

\vspace{1mm}

The above chain of implications is a special case of~\cite[Theorem 5.8]{MR23}, which may be derived from Lemma~\ref{lemma-TPF-MD} as observed in \cite[Remark 6, Remark 17 and Remark 18]{MR23}. As alluded to in Remark \ref{rem-L1andL1NE}, the question whether the equivariant $L^1$-compression exponent is characterized through invariant percolations remains open. More precisely, it is not known that (ii) implies (i) above. This would follow from the ``missing implication'' in \cite[Proposition 2.8]{CSV12}. 

It is also instructive to consider the analogously defined equivariant $L^p$-compression exponent $\alpha_p^{\scriptscriptstyle{\#}}(G)$, $p\in(1,2]$. Towards characterizations, recall from Remark~\ref{rem-Lp} that it is natural to consider the existence of invariant percolations with arbitrarily large marginals and two-point function decay as in~\eqref{equ-LpTPFDecay} with $\rho(t)=t^\alpha$ instead of (ii)~above. With this condition, it follows from similar reasoning as in Remark~\ref{rem-Lp} and as used above, that the chain of implications holds if in (iii) we consider $\alpha_p^{\scriptscriptstyle{\#}}(G)$ and if (i) is replaced with the condition that $\alpha_p^{\scriptscriptstyle{\#}}(G)$ is at least~$2 \alpha$. In other words, the exponent has to be divided by two in both steps. This is due to the fact that we use the square-root trick both to create measure definite kernels from $L^p$-kernels and to create equivariance.

\item[(C)] The present results motivate the investigation of versions of these properties for finite graphs and their local weak limits. Indeed, in percolation theory as well as in coarse geometry, there is a second class of typically non-symmetric graphs of major interest. These are finite graphs and in particular sequences of finite graphs with size tending to infinity, see e.g.~\cite{ABS04,EH24} and \cite{AF19}. These are especially relevant due to their connection with the notion of {\em local weak convergence} introduced in \cite{BS01}, see also \cite{AL07}. The general question about the interplay between percolation and the geometry of such graphs as well as the specific questions investigated in this paper have interesting analogues for finite graphs and their local weak limits. In this context, expanders are well-known examples with poor embedding properties into Hilbert and $L^1$-spaces. In a recent breakthrough development, Salez \cite{S22} has shown using probabilistic arguments that these have negative Ollivier-Ricci curvature. Understanding the relationship between, on the one hand, expansion and Ollivier-Ricci curvature, and, on the other hand, finite and random analogues of the properties considered in this paper could lead to natural refinements of the results in \cite{S22}. We plan to address these questions in future work. 

\item[(D)] In the direction of finding novel links between our methods and random walk approaches (cf.~Remark~\ref{rem-L1ProbInt}), let us mention that the recent preprint \cite{H24} finds an upper bound on the probability that random walk on a group has small total displacement in terms of the spectral and isoperimetric profiles (and provides applications to several open problems related to random walks), the proof of which is inspired by the equivariant method of \cite{MR23}.
\item[(E)] It would be interesting to formulate a non-invariant version of property (T), e.g.~by introducing a weak form of invariance. For Cayley graphs, one possibility could be to consider {\em uniformly quasi-invariant} percolations, i.e.\ general percolations $\P$ such that all shifts admit densities uniformly bounded in $L^\infty(\P)$. The condition in Corollary \ref{cor-CE-groups} with uniformly quasi-invariant percolations would imply that the corresponding group has the {\em weak Haagerup property} introduced in \cite{K14} (this follows from \cite{MR23} together with \cite[Proposition 3.1]{HK15}), but we do not know whether the converse holds. 
\item[(F)] The previous remark seems to be closely related to the challenge of finding a non-trivial, coarse version of property (T) posed by Roe \cite[Section 11.4, Question (b)]{R03}. Note that Proposition~\ref{prop-GeneralDecay} fails for finite graphs. Thus the natural attempt of finding a coarse version of property~(T) by considering the percolation-theoretic condition in (A) for general, instead of invariant, bond percolations is another instance of a plausible definition of "coarse property~(T)" trivializing to "being bounded'', cf.~\cite[Section 11.4.3]{R03}.
\end{itemize}

\medskip

\noindent{\bf Acknowledgement:}
We thank Russell Lyons for sharing insights and very helpful comments on the previous version of the article. The research of both authors is funded by the Deutsche Forschungsgemeinschaft (DFG) under Germany's Excellence Strategy EXC 2044-390685587, Mathematics M\"unster: Dynamics-Geometry-Structure.

\medskip 

\noindent{\bf Funding:} This work was supported by the Deutsche Forschungsgemeinschaft [EXC 2044-390685587].

\noindent{\bf Data Availability:} Data sharing is not applicable to this article as no datasets were generated or analyzed during the current study.

\noindent{\bf Declarations.}

\noindent{\bf Conflict of interest:} The authors declare that they have no conflict of interest.


\begin{thebibliography}{WWW98}

\bibitem{AN84} 
{\sc Aizenman, M.} and {\sc Newman, C. M.} (1984). Tree graph inequalities and critical behavior in percolation models. {\em J. Statist. Phys.} {\bf 36} 107-143.

\bibitem{AL07}
{\sc Aldous, D.} and {\sc Lyons, R.} (2007). Processes on unimodular random networks. {\em Electron. J. Probab.} {\bf 12} 1454-1508.

\bibitem{AF19}
{\sc Alekseev, V.} and {\sc Finn-Sell, M.} (2019). Sofic boundaries of groups and coarse geometry of sofic approximations. {\em Groups Geom. Dyn.} {\bf 13} 191–234.

\bibitem{ABS04}
{\sc Alon, N., Benjamini, I., Stacey, A.} (2004). Percolation on finite graphs and isoperimetric inequalities. {\em Ann. Probab.} {\bf 32} 1727-1745.

\bibitem{ADS09}
{\sc Arzhantseva, G., Drutu, C., Sapir, M.} (2009). Compression functions of uniform embeddings of groups into Hilbert and Banach spaces. {\em J. Reine Angew. Math.} {\bf 633} 213-235.

\bibitem{A11}
{\sc Austin, T.} (2011). Amenable groups with very poor compression into Lebesgue spaces. {\em Duke Math. J.} {\bf159}(2) 187-222.

\bibitem{BdlHV08}
{\sc Bekka, B., de la Harpe, P., Valette, A.} (2008). Kazhdan's property (T). Cambridge University Press, Cambridge.

\bibitem{BLPS99} 
{\sc Benjamini, I., Lyons, R., Peres, Y., Schramm, O.} (1999). Group-invariant percolation on graphs. {\em Geom. Funct. Anal.} \textbf{9} 29-66.

\bibitem{BS01}
{\sc Benjamini, I.} and {\sc Schramm, O.} (2001). Recurrence of distributional limits of finite planar graphs. {\em Electron. J. Probab.} {\bf 6} 1-13.

\bibitem{BZ21}
{\sc Brieussel, J.} and {\sc Zheng, T.} (2021). Speed of random walks, isoperimetry and compression of finitely generated groups. {\em Ann. of Math.} {\bf 193} 1-105.

\bibitem{BO08}
{\sc Brown, N. P.} and {\sc Ozawa, N.} (2008). $C^*$-Algebras and finite dimensional approximations. In: {\em Graduate Studies in Mathematics}, Vol. \textbf{88}, AMS, Providence, Rhode Island.

\bibitem{CDH10}
{\sc Chatterji, I., Drutu, C., Haglund, F.} (2010). Kazhdan and Haagerup properties from the median viewpoint. {\em Adv. Math.} \textbf{225}(2) 882-921.

\bibitem{CMV04}
{\sc Cherix, P.-A., Martin, F., Valette, A.} (2004). Spaces with measured walls, the Haagerup property and property (T). {\em Ergodic Theory Dynam. Sys.} {\bf 24} 1895-1908.

\bibitem{CCJJV01}
{\sc Cherix, P.-A., Cowling, M., Jolissaint, P., Julg, P., Valette, A.} (2001). Groups with the Haagerup Property (Gromov's a-T-menability). In: Progress in Mathematics, Vol. {\bf 197}, Birkhäuser, Boston, MA.

\bibitem{CSV12}
{\sc de Cornulier, Y., Stalder, Y., Valette, A.} (2012). Proper actions of wreath products and generalizations. {\em Trans. Amer. Math. Soc.} {\bf 364}(6) 3159-3184.

\bibitem{DL97}
{\sc Deza, M.} and {\sc Laurent, M.} (1997). Geometry of Cuts and Metrics. Springer, Berlin, 1997.

\bibitem{EH24}
{\sc Easo, P.} and {\sc Hutchcroft, T.} (2024). The critical percolation probability is local. {\em arxiv-preprint:} \url{https://arxiv.org/abs/2310.10983}.

\bibitem{E21b}
{\sc Elek, G.} (2021). Uniform hyperfiniteness. {\em Trans. Amer. Math. Soc.} {\bf 374} 5095-5111.

\bibitem{GK04}
{\sc Guentner, E.} and {\sc Kaminker, J.} (2004). Exactness and uniform embeddability of discrete groups. {\em J. Lond. Math. Soc.} {\bf 70} 703-718.

\bibitem{HK15}
{\sc Haagerup, U.} and {\sc Knudby, S.} (2015). The weak Haagerup property II: Examples. {\em Int. Math. Res. Not. IMRN} {\bf 2015}(16) 6941-6967. 

\bibitem{HJ06}
{\sc H{\"a}ggstr{\"o}m, O.} and {\sc Jonasson, J.} (2006). Uniqueness and non-uniqueness in percolation theory. {\em Probab. Surv.} {\bf 3} 289-344.

\bibitem{H24}
{\sc Hutchcroft, T.} (2024). Small-ball estimates for random walks on groups. {\em arxiv-preprint:} \url{https://arxiv.org/abs/2406.17587}.

\bibitem{K14} 
{\sc Knudby, S.} (2014). Semigroups of Herz-Schur multipliers. {\em J. Funct. Anal.} {\bf 266} 1565-1610.

\bibitem{LP17}
{\sc Last, G.} and {\sc Penrose, M.} (2017). Lectures on the Poisson Process. In: {\em Institute of Mathematical Statistics Textbooks}, Cambridge University Press, Cambridge.

\bibitem{L13}
{\sc Lyons, R.} (2013). Distance covariance in metric spaces. {\em Ann. Probab.} {\bf 41} 3284-3305.

\bibitem{LP16}
{\sc Lyons, R.} and {\sc Peres, Y.} (2016). Probability on Trees and Networks. Cambridge University Press, New York.

\bibitem{LS99} 
{\sc Lyons, R.} and {\sc Schramm, O.} (1999). Indistinguishability of percolation clusters. {\em Ann. Probab.} {\bf 27} 1809-1836.

\bibitem{MR22}
{\sc Mukherjee, C.} and {\sc Recke, K.} (2024). Schur multipliers of $C^*$-algebras, group-invariant compactification and applications to amenability and percolation. {\em J. Funct. Anal.} {\bf 287}(2), 110468.

\bibitem{MR23}
{\sc Mukherjee, C.} and {\sc Recke, K.} (2026). Haagerup property and group-invariant percolation. {\em Probab. Theory Related Fields} {\bf 196} 849–909.

\bibitem{NP08}
{\sc Naor, A.} and {\sc Peres, Y.} (2008). Embeddings of discrete groups and the speed of random walks. {\em Int. Math. Res. Notices IMRN}, rnn076.

\bibitem{NP11}
{\sc Naor, A.} and {\sc Peres, Y.} (2011). $L_p$ compression, traveling salesmen, and stable walks. {\em Duke Math. J.} {\bf 157} 53-108.

\bibitem{P22}
{\sc Pete, G.} (2022). Probability and Geometry on Groups. {\em Available at:} \url{http://math.bme.hu/~gabor/PGG.pdf}.

\bibitem{RS98}
{\sc Robertson, G.} and {\sc Steger, T.} (1998). Negative definite kernels and a dynamical characterization of property (T) for countable groups.  {\em Ergodic Theory Dynam. Sys.} {\bf 18} 247-253.

\bibitem{R03}
{\sc Roe, J.} (2003). Lectures on coarse geometry. In: {\em University Lecture Series}, Vol. \textbf{31}, AMS, Providence, Rhode Island.

\bibitem{S22}
{\sc Salez, J.} (2022). Sparse expanders have negative curvature. {\em Geom. Funct. Anal.} {\bf 32} 1486-1513.

\bibitem{W09}
{\sc Willett, R.} (2009). Some notes on property A. In: {\em Limits of graphs in group theory and computer science}, EPFL Press, Lausanne.

\bibitem{Y00}
{\sc Yu, G.} (2000). The coarse Baum–Connes conjecture for spaces which admit a uniform embedding into Hilbert space. {\em Invent. Math.} {\bf 139} 201-240.

\end{thebibliography}
\end{document}